\documentclass[12pt]{article}

\usepackage{times}
\usepackage[centertags]{amsmath}
\usepackage{amssymb}
\usepackage{amsthm}
\usepackage[latin1]{inputenc}
\usepackage[colorlinks=true,urlcolor=red,citecolor=red,linkcolor=red]{hyperref}
\usepackage[a4paper, top=5cm, bottom=5cm, left=3.5cm, right=2.5cm]{geometry}

\newenvironment{acknowledgements}%
    {\cleardoublepage\thispagestyle{empty}\null\vfill\begin{center}%
    \bfseries Acknowledgements\end{center}}%
    {\vfill\null}

\newcounter{lastnote}

\title{On the codimension of subalgebras of the algebra of matrices over a field}

\author
{Giuseppe Zito\\
\\
\\
}

\date{}

\begin{document}

\baselineskip24pt

\maketitle

\begin{abstract}
In this paper we provide an elementary and easy proof that a proper subalgebra of the matrix algebra $ \mathbb{K}^{n,n }$,  with $n \geq3$  and $\mathbb{K}$ an arbitrary field, has dimension strictly less than $n^2-1$.
\end{abstract}

\section*{Introduction}

 The aim of this  paper is to give   an elementary proof of the following result:
\newtheorem{teorema}{Theorem}
\begin{teorema}
\label{T}
If $\mathcal{A}$ is a proper subalgebra of $ \mathbb{K}^{n,n}$, the algebra of  $ n\times n$ matrices with entries in an arbitrary field $\mathbb{K}$, and if $n\geq 3$, then  $\dim(\mathcal{A}) \leq n^2-2$, where $\dim(\mathcal{A})$ is the dimension of $\mathcal{A}$ as a vector space over $\mathbb{K}$.
\end{teorema} 
When $\mathbb{K}$ is an algebraically closed field, this result can be deduced and  generalized by the following well known Theorem of Burnside (see \cite{HR}, \cite{Lam}, \cite{LR} and \cite{Ros} for elementary proofs):
\newtheorem*{teorema2}{Burnside's Theorem}
\begin{teorema2}
If $\mathcal{A}$ is an irreducible subalgebra of $ \mathbb{K}^{n,n}$, with $n\geq 1$ and $\mathbb{K}$ is an algebraically closed field, then $\mathcal{A}=\mathbb{K}^{n,n}$. 
\end{teorema2} 

Recall that a subalgebra $\mathcal{A} \subseteq \mathbb{K}^{n,n}$ is said to be irreducible if the only linear subspaces $U\subseteq \mathbb{K}^n$ invariant under all the elements of $\mathcal{A}$, i.e. such that $AU \subseteq U$ for all $A \in \mathcal{A}$, are $\left\{0 \right\}$ and $\mathbb K^n$.

Then, if $\mathcal{A}$ is a proper subalgebra of $\mathbb{K}^{n,n}$ , with  $\mathbb{K}$ an algebraically closed field, there exists a non trivial subspace $U$ of $\mathbb{K}^n$ 
that is invariant under all the elements of $\mathcal{A}$ with $0<\dim(U)=r<n$. This immediately implies $\dim(\mathcal A)\leq n^2-r(n-r)\leq n^2-n+1$, refining significantly the bound
in Theorem \ref{T}.



So  Theorem \ref{T}   gives less information on the dimension of a proper subalgebra of $\mathbb{K}^{n,n}$ than the consequences of Burnside's Theorem, but it holds over an arbitrary field.

From another perspective the reason for which there do not exist codimension 1 subalgebras of $\mathbb K^{n,n}$ resides on Wedderburn's Theorem providing the decomposition of $\mathbb K^{n,n}$ into the direct
sum of its {\it semi-simple} and {\it radical} parts. Over an algebraically closed field the semi simple part is a direct sum of matrix algebras while the dimension of the radical is controlled
by Gerstenhaber's Theorem on the dimension of linear subspaces of nilpotent matrices, see \cite{Ger}. Putting together  these not  trivial facts one obtains
a different proof of the above inequality $\dim(\mathcal A)\leq n^2-n+1$.  This approach is present in \cite{Agore}, where this inequality is also proved to be true  over a field  of characteristic zero, thus requiring less restrictive hypothesis about $\mathbb K$. A, not trivial, proof  of Theorem \ref{T} with $\mathcal A$ over an arbitrary field, can be also found as a consequence of the results contained in \cite{Militaru}.

 In conclusion we post our novelty only on the simple and elementary proof of Theorem \ref{T} but not on its contents which are surely
well known to any expert in the field.

\section*{Proof of Theorem \ref{T}}

Suppose that $\dim(\mathcal A)=n^2-1$.  Since $ \dim(\mathbb{K}^{n,n})=n^2$ ,  given a basis for the vector space  $ \mathbb{K}^{n,n}$, the subspace $\mathcal A$ is represented by only one  homogeneous equation in the associated coordinates. 
We consider the standard basis  of $ \mathbb{K}^{n,n}$ which consists of the matrices $E_{ij}$ such that $\left( E_{ij} \right)_{k,l}=\delta_{ik}\cdot \delta_{jl}$, where $\delta$ is the Kronecker delta.

We have:
$$ \left(E_{ij} \cdot E_{kl} \right)_{m,q}=\sum_{p=1}^n{ \left( E_{ij}\right)_{m,p}\left( E_{kl}\right)_{p,q}}=\sum_{p=1}^n{\delta_{im}\delta_{jp}\delta_{kp}\delta_{lq}}=\delta_{im}\delta_{kj}\delta_{lq}.$$
If $k \neq j$ then $E_ {ij} \cdot E_{kl}=0$, while, if $k=j$, from $$ \left(E_{ij} \cdot E_{jl} \right)_{m,q}=\delta_{im} \delta_{lq}= \left( E_{il} \right)_{m,q},$$ 
it follows that  $E_ {ij} \cdot E_{jl}=E_{il}.$ 

We now suppose the subspace $\mathcal A$ is  given by the following equation:
\begin{equation}  \label{equazione} a_{11}x_{11}+ \ldots + a_{1n}x_{1n}+a_{21}x_{21}+\ldots + a_{nn}x_{nn}=0 \end{equation}
where $a_{ij} \in \mathbb{K}$ for all $i,j=1,\ldots,n$.

We want to show that $a_{ij}=0$ for all $i$ and $j$. To this aim, for all $i, j, k$ and $l \in \left\{1, \ldots, n \right\}$, we define the matrices $$ D_{ij}^{kl}=a_{kl}E_{ij}-a_{ij}E_{kl}.$$
By construction $D_{ij}^{kl}$ satisfies the equation of $\mathcal A$, and therefore we have $D_{ij}^{kl} \in \mathcal{A}$ for all $i, j, k$ and $l \in \left\{1, \ldots, n \right\}$ .

We now consider a triple of pairwise distinct indices $i$, $j$ and $k$ (we can find them because $n \geq 3$).  $D_{ik}^{ij}$ and $D_{kj}^{ij}$ belong to $\mathcal{A}$, which is closed under  matrix multiplication.  Therefore we  deduce:

$$ \mathcal{A} \ni D_{ik}^{ij} \cdot D_{kj}^{ij}=\left( a_{ij}E_{ik}-a_{ik}E_{ij}\right) \cdot \left( a_{ij}E_{kj}-a_{kj}E_{ij}\right)=$$ $$=a_{ij}^2 E_{ik}E_{kj}-a_{ij}a_{kj}E_{ik}E_{ij}-a_{ik}a_{ij}E_{ij}E_{kj}+a_{ik}a_{kj}E_{ij}^2=a_{ij}^2 E_{ij},$$

where the last equality follows from the properties of the matrices $E_{ij}$ and from the choice of $i,j$ e $k$.
The matrix $a_{ij}^2 E_{ij}$ belongs to $\mathcal A$, so it must satisfy the equation (\ref{equazione}).
From this it follows that $a_{ij}^3=0 \Rightarrow a_{ij}=0$.

We have just proved that $a_{ij}=0$ for all $i\neq j$, so the equation of $\mathcal A$ can be written in the form:

\begin{equation} \label{equazione due} a_{11}x_{11}+a_{22}x_{22}+\ldots +a_{nn}x_{nn}=0. \end{equation}

From this equation we can deduce that, for all $i=1, \ldots,n$ and $k\neq i$, the matrices $P_i=a_{ii}E_{ik}$ and $Q_i=a_{ii} E_{ki}$ belong to $\mathcal A$. 

Then we have:
$$ P_i \cdot Q_i= a_{ii}^2 E_{ik} E_{ki}= a_{ii}^2 E_{ii} \in \mathcal{A}. $$

The matrix $a_{ii}^2 E_{ii}$ belongs to $\mathcal A$, so it must satisfy the equation (\ref{equazione due}).
From this it follows that $a_{ii}^3=0 \Rightarrow a_{ii}=0$ for all $i=1, \ldots n$.

We proved that  $a_{ij}=0$ for all $i$ and $j$, in contradiction with the assumption that $\mathcal A$ can be represented by only one equation.

        \begin{acknowledgements}
    I would like to thank the editorial board of  " Esercitazioni Matematiche ", an initiative of the Department of Mathematics and Computer Science of Catania, which introduced me to this problem,  encouraging  me to solve it. 

In particular I would like to thank Professor Francesco Russo for his constructive comments during the realization of this paper.
        \end{acknowledgements}

E-mail address: giuseppezito@hotmail.it


\begin{thebibliography}{Dillo 100}

\bibitem{Agore}  A.~L.~Agore,  {\em The maximal dimension of unital subalgebras of the matrix algebra}, arXiv:14030773.
\bibitem{Militaru} A.~L.~Agore, G.~Militaru, {\em  Extending structures, Galois groups and supersolvable associative algebras}, accepted for publication in Monatshefte fur Mathematik.

\bibitem{Ger} M.~Gerstenhaber,  {\em
On nilalgebras and linear varieties of nilpotent matrices}, in I. Amer. J. Math,  {\em 80} (1958),  614-622.

\bibitem{HR} I.~Halperin, P.~Rosenthal,  {\em
Burnside's theorem on algebras of matrices}, in Amer. Math. Monthly,  {\em 87} (1980), 810. 
\bibitem{Lam} T.~Y.~Lam,  {\em
A theorem of Burnside on matrix rings}, in Amer. Math. Monthly,  {\em 105} (1998),  651-653.

\bibitem{LR} V.~ Lomonosov,  P.~Rosenthal,  {\em
The simplest proof of Burnside's theorem on matrix algebras}, in Linear Algebra  Appl., {\em 383} (2004),  45-47.


\bibitem{Ros} E.~Rosenthal,  {\em
A remark on Burnside's theorem on matrix algebras}, in Linear Algebra  Appl.,  {\em 63} (1984), 175-177.


\end{thebibliography}
\end{document}